\documentstyle[11pt]{article}\input{amssym.def}\input{amssym}
\begin{document}
\def\b1{\text{\bf 1}}\def\bF{\text{\bf F}}\def\br{\text{\bf r}}\def\bv{\text{\bf v}}
\def\BA{{\Bbb A}}\def\BC{{\Bbb C}}\def\BN{{\Bbb N}}\def\BP{{\Bbb P}}\def\BR{{\Bbb R}}
\def\BZ{{\Bbb Z}}\def\CA{{\cal A}}\def\CB{{\cal B}}\def\CD{{\cal D}}\def\CDiff{{\cal{D}}iff}
\def\CE{{\cal E}}\def\CF{{\cal F}}\def\CG{{\cal G}}\def\CHom{{\cal{H}}om}\def\CL{{\cal L}}
\def\CM{{\cal M}}\def\CN{{\cal N}}\def\CO{{\cal O}}\def\CP{{\cal P}}\def\CR{{\cal R}}
\def\CT{{\cal T}}\def\CU{{\cal U}}\def\CV{{\cal V}}\def\CW{{\cal{W}}}
\def\dq{\buildrel{.}\over{q}}\def\dbr{\buildrel{.}\over{\text{\bf r}}}
\def\dbv{\buildrel{.}\over{\text{\bf v}}}\def\ddbr{\buildrel{..}\over{\text{\bf r}}}
\def\Der{\text{Der}}\def\dpar{\partial}\def\End{\text{End}}
\def\fB{{\frak B}}\def\fC{{\frak C}}\def\fc{{\frak c}}\def\fD{{\frak D}}\def\ff{{\frak f}}
\def\fg{{\frak g}}\def\fG{{\frak G}}\def\fq{{\frak q}}\def\fR{{\frak R}}\def\fk{{\frak k}}
\def\ft{{\frak t}}\def\gl{{\text{gl}}}\def\Gr{{\text{Gr}}}
\def\hA{\widehat A}\def\hCO{\widehat{{\cal O}}}\def\hCT{\widehat{{\cal T}}}\def\hD{{\widehat D}}
\def\hfg{\widehat{\frak g}}\def\hgamma{\widehat{\gamma}}\def\hgl{\widehat{\text{gl}}}
\def\hO{{\widehat O}} \def\hOmega{\widehat{\Omega}}\def\hsl{\widehat{\text{sl}}}
\def\htau{\widehat{\tau}}\def\hV{\widehat V}\def\hX{\widehat X}\def\Id{\text{Id}}\def\id{\text{id}}
\def\Ker{\text{Ker}}\def\Lie{\text{Lie}}\def\MQ{{\cal{MQ}}}\def\QM{{\cal{QM}}}\def\Res{\text{Res}}
\def\Spec{\text{Spec}}\def\Spf{\text{Spf}}\def\ta{\tilde a}\def\tb{\tilde b}\def\tc{\tilde c}
\def\tCT{\widetilde{\cal{T}}}\def\td{\tilde d}\def\tF{\tilde{F}}\def\tf{\tilde f}\def\tG{\tilde G}
\def\tg{\tilde g}\def\tI{\tilde I}\def\tJ{\tilde J}\def\tK{\tilde K}\def\tL{\tilde L}
\def\tphi{\tilde\phi}\def\tpsi{\tilde\psi}\def\tQ{\tilde Q}\def\Tr{\text{Tr}}\def\tS{\tilde{S}}
\def\ttau{\tilde{\tau}}\def\tW{\tilde{W}}\def\tX{\tilde{X}}\def\tx{\tilde x}\def\txi{\tilde\xi}
\def\ty{\tilde y}\def\tz{\tilde z}\def\Vir{{\cal V}ir}
\def\btu{\bigtriangleup}\def\hra{\hookrightarrow}\def\iso{\buildrel\sim\over\longrightarrow} 
\def\lla{\longleftarrow}\def\lra{\longrightarrow}\def\ra{\rightarrow}
\def\twolra{\buildrel\longrightarrow\over\longrightarrow}\def\deg{\text{deg}}
\def\deg{{\mathrm{deg}}\,}\def\binom#1#2{\left(\!\!\begin{array}{c} #1\\ #2\end{array}\!\!\right)}

\hfill CPT-P09-2006

\bigskip\bigskip\bigskip

\centerline{\Large UNE INTERSECTION DE QUADRIQUES LI\'EE \`A LA SUITE DE STURM}

\bigskip\bigskip

\centerline{Oleg Ogievetsky\footnote{Centre de Physique Th\'eorique, Luminy, 13288 Marseille, 
France (Unit\'e Mixte de Recherche 6207 du CNRS et des Universit\'es Aix--Marseille I, 
Aix--Marseille II et du Sud Toulon -- Var; laboratoire affili\'e \`a la FRUMAM, FR 2291) et  
Institut de Physique P.N. Lebedev, D\'epartement Th\'eorique, Leninsky prospekt 53, 119991 Moscou, 
Russie; \hfill \\ email: oleg@cpt.univ-mrs.fr}, Vadim Schechtman\footnote{Laboratoire Emile Picard, 
UFR MIG, Universit\'e Paul Sabatier, 31062 Toulouse, France; \hfill \\
email: schechtman@math.ups-tlse.fr}}

\vskip 3cm

\centerline{TABLE DES MATI\`ERES} 

\bigskip\bigskip

Premi\`ere Partie.\ Formules  

\medskip

\ \ \ \ \S\  1.\ \ \  Introduction\hfill\pageref{Introduction} \ \ \ 

\ \ \ \ \S\  2.\ \ \  Alg\`ebre $\fB$\hfill\pageref{Algebre} \ \ \ 

\ \ \ \ \S\  3.\ \ \  D\'ebut de la d\'emonstration du th\'eor\`eme 1.5\hfill\pageref{Debut} \ \ \ 

\ \ \ \ \S\  4.\ \ \  Formule $(A)$\hfill \pageref{FormuleA} \ \ \ 

\ \ \ \ \S\  5.\ \ \  Formule $(B)$\hfill \pageref{FormuleB} \ \ \ 

\bigskip

Deuxi\`eme Partie.\ Polyn\^omes d'Euler et d\'eterminant de Cauchy

\medskip

\ \ \ \ \S\  1.\ \ \  Nombres $\beta(j)_i$\hfill\pageref{Nombresbeta} \ \ \ 

\ \ \ \ \S\  2.\ \ \  Polyn\^omes d'Euler et fonction 
hyperg\'eom\'etrique\hfill\pageref{PolynomesEuler} \ \ \  
 
\ \ \ \ \S\  3.\ \ \  Asymptotiques\hfill\pageref{Asymptotiques} \ \ \ 

\bigskip

Bibliographie\hfill\pageref{Bibliographie} \ \ \

\newpage 
\centerline{PREMI\`ERE PARTIE. FORMULES}

\bigskip\bigskip 

\centerline{\bf \S\  1. Introduction}\label{Introduction}

\bigskip

{\bf 1.1.} Cet article est une variation sur un th\`eme de [Jacobi].  

Soit 
$$ f(x) = a_n x^n + a_{n-1}x^{n-1} + \ldots + a_0 $$
un polyn\^ome de degr\'e $n > 0$ \`a coefficients dans un corps de 
base $\fk$ de caract\'eristique $0$.  
Rappelons que {\it la suite de Sturm} de $f$,  
$$ \ff = (f_0, f_1, f_2, \ldots)\ , $$ 
est d\'efinie par r\'ecurrence~: on pose $f_0(x) = f(x),\ f_1(x) = f'(x)$ et 
pour $j \geq 1$ $f_{j+1}$ est le reste de la division euclidienne 
de $f_{j-1}$ par $f_{j}$, avec le signe oppos\'e~: 
$$ f_{j-1}(x) = q_{j-1}(x)f_j(x) - f_{j+1}(x),\ \eqno{(1.1.1)}$$
$\deg f_{j+1}(x) < \deg f_{j}(x)$, cf. le c\'el\`ebre m\'emoire [Sturm]. 

Dans cette note on propose des formules explicites pour les coefficients 
des polyn\^omes $f_j$ en termes des coefficients de $f$. 
Plus g\'en\'eralement, on donnera des formules analogues pour les membres de l'algorithme 
d'Euclide correspondant \`a deux polyn\^omes quelconques $f_1, f_2$ de degr\'es 
$n-1, n-2$. 

Notre point de d\'epart est une alg\`ebre $\fB$, quotient de l'anneau de
polyn\^omes en variables $b(i)_j\ (i\geq 1,\ j\geq 2i)$ par certains 
r\'elations quadratiques, cf. (1.7.1) ci-dessous. Nos formules sont des 
cons\'equences 
des identit\'es dans $\fB$, analogues des r\'elations de Pl\"ucker.

\medskip
{\bf 1.2.} Pour \'enoncer le r\'esultat, introduisons les quantit\'es 
quadratiques
$$ b(j)_i = n\sum_{p=0}^{j-1}\ (i-2p)a_{n-p}a_{n-i+p} 
- j(n-i+j)a_{n-j}a_{n+j-i}, $$
$j \geq 1,\ i \geq 2j$. Ici on pose $a_i = 0$ pour $i < 0$. Par exemple,  
$$ b(1)_i = n i a_n a_{n-i} - (n-i+1)a_{n-1}a_{n-i+1}\ . $$

\medskip
{\bf 1.3.} Ensuite on introduit, pour $m\geq 2$, les matrices $(m-1)\times (m-1)$ 
sym\'etriques  
$$ C(m) = \left(\begin{array}{cccccc} 
b(1)_2 & b(1)_3 & b(1)_4 & b(1)_5 & \ldots & b(1)_m \\
b(1)_3 & b(2)_4 & b(2)_5 & b(2)_6 & \ldots & b(2)_{m+1} \\
b(1)_4 & b(2)_5 & b(3)_6 & b(3)_7 & \ldots & b(3)_{m+2} \\
b(1)_5 & b(2)_6 & b(3)_7 & b(4)_8 & \ldots & b(4)_{m+3} \\
. & . & . & . & \ldots & . \\ 
b(1)_m & b(2)_{m+1} & b(3)_{m+2} & b(4)_{m+3} & \ldots & b(m-1)_{2m-2} \\ 
\end{array}\right)\ . $$

De plus, pour $i \geq 0$ on d\'efinit une matrice "d\'ecal\'ee" $C(m)_i$~: 
elle est obtenue en rempla\c{c}ant dans $C(m)$ la derni\`ere ligne par
$$ \left(\begin{array}{cccccc}
b(1)_{m+i} & b(2)_{m+i+1} & b(3)_{m+i+2} & b(4)_{m+i+3} & \ldots & b(m-1)_{2m+i-2} \\ 
\end{array}\right)\ . $$
Donc $C(m)_0 = C(m)$. On pose 
$$ c(m)_i := \det C(m)_i,\ c(m) := c(m)_0\ . $$
En particulier, 
$$ c(2)_i = b(1)_{i+2} $$

Il est commode de poser 
$$ c(1)_i := \frac{(n-i)a_{n-i}}{na_n}, $$
$i \geq 0$, $c(1):= c(1)_0 = 1$.  

\medskip
{\bf 1.4.} Puis on d\'efinit les nombres $\gamma_j,\ j \geq 1$ par 
r\'ecurrence~:
$$ \gamma_1 = na_n,\ \gamma_2 = - \frac{1}{n^2a_n},\ 
\gamma_{j+1} = \gamma_{j-1}\cdot \frac{c(j-1)^2}{c(j)^2}, $$
$j \geq 2$. Autrement dit, 
$$ \gamma_j = (-1)^{j+1}\epsilon_j\cdot \prod_{i=1}^{j-2}\ c(j-i)^{2(-1)^i} \ ,$$
o\`u $\epsilon_j = n a_n$ si $j$ est impair et $1/(n^2a_n)$ sinon. 

Les nombres  $\gamma_1,\ldots,\gamma_j$ sont donc bien d\'efinis si tous 
les nombres $c(2), c(3), \ldots, c(j-1)$ sont diff\'erents de z\'ero. 

\medskip
{\bf 1.5.} {\it Th\'eor\`eme.} Supposons que $\deg f_j = n-j$, donc 
$\deg f_i = n - i$ pour $i \leq j$. 

Alors pour tous $i \leq j$, on a $c(i) \neq 0$ et 
$$ f_i(x) = \gamma_i \cdot \sum_{p=0}^{n-i}\ c(i)_px^{n-i-p} \ . $$
En particulier, le coefficient dominant de $f_i(x)$ est \'egal \`a $\gamma_i c(i)$. 

\medskip
{\bf 1.6.} On v\'erifie aussit\^ot que 
$$ b(k)_i - b(k-1)_i = c(1)_{k-1}b(1)_{i-k+1} - c(1)_{i-k}b(1)_k
\eqno{(1.6.1)} $$
pour tous $k \geq 2,\ i \geq 2k-2$. Par exemple, 
$$ b(2)_i - b(1)_i = c(1)_{1}b(1)_{i-1} - c(1)_{i-2}b(1)_2, $$   
$$ b(3)_i - b(2)_i = c(1)_{2}b(1)_{i-2} - c(1)_{i-3}b(1)_3, $$
etc. Il s'en suit que tous les $b(j)_i,\ j \geq 2,$ sont expressibles en termes de 
$c(1)_p$ et $c(2)_p = b(1)_{p+2},\ p \geq 0$. 

\medskip
{\bf 1.7.} Les formules (1.6.1) impliquent que 
les nombres $b(i)_j$ satisfont aux relations quadratiques suivantes~: 
$$\begin{array}{c} \bigl(b(k)_i - b(k-1)_i\bigr)\cdot b(1)_j \\[1em]
 = \bigl(b(j)_{i-k+j} - b(j-1)_{i-k+j}\bigr)\cdot b(1)_k 
- \bigl(b(j)_{k+j-1} - b(j-1)_{k+j-1}\bigr)\cdot b(1)_{i-k+1}
\end{array}\eqno{(1.7.1)} $$
On verra que la preuve de 1.5 ne d\'epend que des relations (1.7.1). 
 
On formalise la situation en introduisant une alg\`ebre quadratique 
correspondante, cf. \S\  2 ci-dessous. 

\medskip
{\bf 1.8.} Maintenant soient 
$$ f_1(x) = \alpha_0x^{n-1} + \alpha_1x^{n-2} + \ldots $$
et 
$$ f_2(x) = \beta_0x^{n-2} + \beta_1x^{n-3} + \ldots $$
deux polyn\^omes arbitraires de degr\'es $n-1, n-2$. On d\'efinit 
$f_j,\ j \geq 3$ \`a partir de $f_1, f_2$ par les formules 
de l'algorithme d'Euclide (1.1.1).  

Posons
$$ 
c(1)_i := \frac{\alpha_i}{\alpha_0},\ b(1)_{i+2}:= \beta_i,\ i \geq 0\ .
$$
{\it D\'efinissons} les nombres $b(k)_i,\ k \geq 2$ par r\'ecurrence sur 
$k$, \`a partir des formules (1.6.1). 

D\'efinissons les nombres $c(m)_i,\ m\geq 2,$ par les formules 1.3. 

Enfin, on pose~:        
$$
\tilde{\gamma}_1=\alpha_0\ ,\  \tilde{\gamma}_2=1\ ,\ 
\tilde{\gamma}_{j+1}= \tilde{\gamma}_{j-1}\frac{c(j-1)^2}{c(j)^2}\  
$$

Alors on a 

\medskip
{\bf 1.9.} {\it Th\'eor\`eme.} Supposons que $\deg f_j = n-j$, d'o\`u  
$\deg f_i = n - i$ pour $i \leq j$. 

Alors pour tous $i \leq j$, on a $c(i) \neq 0$ et 
$$ f_i(x) = \tilde\gamma_i \cdot \sum_{p=0}^{n-i}\ c(i)_px^{n-i-p} \ .$$
En particulier, le coefficient dominant de $f_i(x)$ est \'egal \`a 
$\tilde\gamma_i c(i)$. 

Cf. [Jacobi], section 15.  

\medskip
{\bf 1.10.} Dans la Deuxi\`eme Partie on pr\'esente un exemple num\'erique. 
L\`a, les d\'eterminants de Cauchy apparaissent dans les asymptotiques des 
coefficients dominants de la suite de Sturm pour les polyn\^omes d'Euler.  

\bigskip\bigskip
%\newpage

\centerline{\bf \S\  2. Alg\`ebre $\fB$}\label{Algebre} 

\bigskip 

{\bf 2.1.} On peut r\'e\'ecrire les relations (1.7.1) sous la forme suivante~: 
$$\begin{array}{l} \det\left(\begin{array}{cc} b(1)_j & b(1)_k \\ 
b(j-1)_{i+j-k} & b(k-1)_i \\ \end{array}\right) - 
\det\left(\begin{array}{cc} b(1)_j & b(j-1)_{j+k-1} \\ 
b(1)_{i-k+1} & b(k)_i \\ \end{array}\right)\\[1em]\ \ \ \ \ 
 + \det\left(\begin{array}{cc} b(1)_k & b(j)_{j+k-1} \\ 
b(1)_{i-k+1} & b(j)_{i+j-k} \\ \end{array}\right) 
 = \Delta(k,j)_i - \Delta'(k,j)_i + \Delta''(k,j)_i = 0\ . \end{array}\eqno{(2.1.1)} $$

\medskip
{\bf 2.2.} On d\'efinit une alg\`ebre quadratique $\fB$ comme une 
$\fk$-alg\`ebre commutative engendr\'ee par les lettres $b(i)_j,\ 
i,j\in \BZ$, modulo les relations (2.1.1), o\`u $i,j,k \in \BZ$. 

(D'ailleurs, dans tout le paragraphe qui suit on peut 
remplacer le corps de base $\fk$ par un anneau commutatif quelconque.) 

\medskip
{\bf 2.3.} Le but de ce paragraphe est d'\'ecrire certaines 
relations entre les d\'eterminants 
$n \times n$ dans $\fB$ qui g\'en\'eralisent (2.1.1). 

On fixe  un nombre entier $n \geq 2$.  Soient $m_1,\ldots,m_{n}, i$  
des entiers. 

On d\'efinit $2n+2$ vecteurs $v_j, w_j \in \fk^{n},\ j = 1, \ldots, n+1$~:      
$$ w_1 = (b(1)_{m_1}, b(1)_{m_2}, \ldots , b(1)_{m_{n}}) \ ,$$
$$ w_{j+1} = (b(1)_{m_1},  \ldots , \hat b(1)_{m_{n+1-j}},\ldots, 
b(1)_{m_{n}}, b(1)_{i - m_{n} + 1}), $$
$1 \leq j \leq n$   
(suivant l'usage, $\hat x$ signifie que l'on omet la composante $x$). 

Puis
$$ v_1 = (b(m_1 - 1)_{i+m_1-m_n}, b(m_2 - 1)_{i+m_2-m_n}, \ldots, 
b(m_{n-1} - 1)_{i + m_{n-1} - m_{n}}, b(m_n-1)_i) \ ,$$ 
$$ v_2 = (b(m_1 - 1)_{m_1+m_n-1}, b(m_2 - 1)_{m_2+m_n-1}, \ldots, 
b(m_{n-1} - 1)_{m_{n-1}+m_{n}-1}, b(m_n)_i) \ ,$$ 
$$ v_3 = (b(m_1 - 1)_{m_1+m_{n-1}-1}, b(m_2 - 1)_{m_2+m_{n-1}-1}, \ldots, 
b(m_{n-2} - 1)_{m_{n-2}+m_{n-1}-1},  $$
$$ b(m_{n-1})_{m_{n-1}+m_{n}-1},b(m_{n-1})_{i+m_{n-1}-m_{n}} )\ ,$$ 
$$ v_4 = (b(m_1 - 1)_{m_1+m_{n-2}-1}, b(m_2 - 1)_{m_2+m_{n-2}-1}, \ldots, 
b(m_{n-3} - 1)_{m_{n-3}+m_{n-2}-1}, $$
$$ b(m_{n-2})_{m_{n-2}+m_{n-1}-1},b(m_{n-2})_{m_{n-2}+m_{n}-1}, 
b(m_{n-2})_{i+m_{n-2} - m_n} )\ , $$
$$ .\ .\ .\ $$
$$ v_n = (b(m_1 - 1)_{m_1+m_{2}-1},b(m_2)_{m_2+m_3-1}, 
b(m_2)_{m_2+m_4-1},\ldots, b(m_2)_{m_2+m_n-1}, b(m_2)_{i+m_2-m_n}) $$
$$ v_{n+1} = (b(m_1)_{m_1+m_{2}-1},b(m_1)_{m_1+m_3-1}, 
\ldots, b(m_1)_{m_1+m_n-1}, b(m_1)_{i+m_1-m_n})\ . $$

\medskip
{\bf 2.4.} Soit 
$$ M = \left(\begin{array}{ccc} x_{11} & \ldots & x_{1,n+1} \\ 
. & \ldots & . \\
x_{n-2,1} & \ldots & x_{n-2,n+1} \\ 
\end{array}\right) $$
une matrice $(n-2)\times (n+1)$ sur $\fB$~; soit $M_i,\ i = 1,\ldots, n+1$, 
ses sous-matrices $(n-2)\times n$. Pour \'ecrire $M_i$, on enl\`eve donc 
la $i$-i\`eme colonne de $M$. 

Maintenant on va d\'efinir $n+1$ matrices $n\times n$ 
$$ D_j = D_j(m_1,\ldots,m_n;M_{n+2-j})_i, $$
$j = 1,\ldots,n+1$. On pose~: 
$$ D_1 = \left(\begin{array}{c} w_1 \\ 
M_{n+1} \\ v_1 \\ \end{array}\right), 
D_j = \left(\begin{array}{ccc} w_j^t & M_{n+2-j}^t & v_j^t \\ 
\end{array}\right),  $$
$j = 2, \ldots, n+1$. Ici $(.)^t$ d\'esigne la matrice transpos\'ee.   

Enfin, on pose 
$$ \Delta_j = \Delta_j(m_1,\ldots,m_n;M_{n+2-j})_i = 
\det D_j(m_1,\ldots,m_n;M_{n+2-j})_i, $$
$j = 1,\ldots,n+1$. 

Consid\'erons la somme altern\'ee 
$$ R(n;m_1,\ldots,m_n;M)_i = \sum_{j=1}^{n+1}\ (-1)^{j+1} 
\Delta_j(m_1,\ldots,m_n;M_{n+2-j})_i\ .$$

\medskip
{\bf 2.5.} {\it Exemple.} $n=2$. Dans ce cas il n'y a pas de matrice $M$~; 
trois nombres entiers sont donn\'es~: $m_1, m_2$ et $i$. On aura $6$ vecteurs~: 
$$ w_1 = (b(1)_{m_1}, b(1)_{m_2}),\ 
w_2 = (b(1)_{m_1}, b(1)_{i-m_2+1}),\ 
w_3 = (b(1)_{m_2}, b(1)_{i-m_2+1}) $$
et
$$ v_1 = (b(m_1 - 1)_{i+m_1-m_2}, b(m_2-1)_i),\  
v_2 = (b(m_1 - 1)_{m_1+m_2-1}, b(m_2)_i), $$
$$ v_2 = (b(m_1)_{m_1+m_2-1}, b(m_1)_{i+m_1-m_2}) \ .$$
Il s'ensuit~:
$$ R(2;m_1,m_2)_i =  \det\left(\begin{array}{cc} b(1)_{m_1} & b(1)_{m_2} \\ 
b(m_1-1)_{i+m_1-m_2} & b(m_2-1)_i \\ 
\end{array}\right) $$
$$ - \det\left(\begin{array}{cc} b(1)_{m_1} & b(m_1-1)_{m_1+m_2-1} \\ 
b(1)_{i-m_2+1} & b(m_2)_i \\   
\end{array}\right) + 
\det\left(\begin{array}{cc} b(1)_{m_2} & b(m_1)_{m_1+m_2-1} \\ 
b(1)_{i-m_2+1} & b(m_2)_{i+m_1-m_2} \\
\end{array}\right) $$
On reconna\^\i t l\`a la partie gauche de (2.1.1) pour $(j,k) = (m_1,m_2)$. 
Il en d\'ecoule que $R(2;m_1,m_2)_i = 0$.  

\medskip
{\bf 2.6.} {\it Exemple.} $n=3$. Dans ce cas la matrice $M$ se r\'eduit \`a 
$4$ \'el\'ements~: 
$$ M = \left(\begin{array}{cccc} x_1 & x_2 & x_3 & x_4 \\ \end{array}\right)\ .$$
L'expression $R(3;m_1,m_2,m_3;M)_i$ prend la forme
$$\begin{array}{l} R(3;m_1,m_2,m_3;M)_i = \det\left(\begin{array} {ccc}
b(1)_{m_1} & b(1)_{m_2} & b(1)_{m_3} \\ 
x_1 & x_2 & x_3 \\ 
b(m_1-1)_{i+m_1-m_3} & b(m_2-1)_{i+m_2-m_3} & b(m_3-1)_i \\ 
\end{array}\right) \\[2em]\ \ \ \ \ 
- \det\left(\begin{array}{ccc} b(1)_{m_1} & x_1 & b(m_1-1)_{m_1+m_3-1} \\ 
b(1)_{m_2} & x_2 & b(m_2-1)_{m_2+m_3-1} \\ 
b(1)_{i-m_3+1} & x_4 & b(m_3)_{i} \\
\end{array}\right) + 
\det\left(\begin{array}{ccc} b(1)_{m_1} & x_1 & b(m_1-1)_{m_1+m_2-1} \\
b(1)_{m_3} & x_3 & b(m_2)_{m_2+m_3-1} \\ 
b(1)_{i-m_3+1} & x_4 & b(m_2)_{i+m_2-m_3} \\
\end{array}\right) \\[2em]\ \ \ \ \ 
- \det\left(\begin{array}{ccc} b(1)_{m_2} & x_2 & b(m_1)_{m_1+m_2-1} \\
b(1)_{m_3} & x_3 & b(m_1)_{m_1+m_3-1} \\ 
b(1)_{i-m_3+1} & x_4 & b(m_1)_{i+m_1-m_3} \\
\end{array}\right)\ .\end{array}$$
Calculons cette expression. 

On d\'eveloppe le premier d\'eterminant suivant la deuxi\`eme 
ligne et les autres suivant les deuxi\`emes colonnes~: 
$$ \Delta_1(3;m_1,m_2,m_3;M_4)_i = 
- x_1 \Delta_1(2;m_2,m_3)_i + x_2 \Delta_1(2;m_1,m_3)_i  
- x_3 \Delta_1(2;m_1,m_2)_{i+m_2-m_3}\ ,$$
$$ \Delta_2(3;m_1,m_2,m_3;M_3)_i = 
- x_1 \Delta_2(2;m_2,m_3)_i + x_2 \Delta_2(2;m_1,m_3)_i  
- x_4 \Delta_1(2;m_1,m_3)_{m_2+m_3-1}\ .$$
Puis
$$ \Delta_3(3;m_1,m_2,m_3;M_2)_i = 
- x_1 \Delta_3(2;m_2,m_3)_i + x_3 \Delta_2(2;m_1,m_2)_{i+m_2-m_3}  
- x_4 \Delta_2(2;m_1,m_3)_{m_2+m_3-1}  $$
et 
$$\Delta_4(3;m_1,m_2,m_3;M_1)_i = 
- x_2 \Delta_3(2;m_1,m_3)_i + x_3 \Delta_3(2;m_1,m_2)_{i+m_2-m_3}  
- x_4 \Delta_3(2;m_1,m_3)_{m_2+m_3-1}\ .$$
Pour abr\'eger les notations on introduit des vecteurs entiers~:
$$(i_1, i_2, i_3, i_4) := (i, i, i + m_2 - m_3, m_2 + m_3 - 1)\ ,$$ 
$$\mu = (m_1,m_2,m_3)\ ,$$
$$\mu_1 = (m_2,m_3),\ \mu_2 = (m_1,m_3),\ \mu_3 = (m_1,m_2)\ .$$  

On peut r\'e\'ecire les formules ci-desssus sous une forme matricielle~: 
$$\left(\begin{array}{c} \Delta_1(3;\mu;M_4)_i \\ 
- \Delta_2(3;\mu;M_3)_i \\ 
\Delta_3(3;\mu;M_2)_i \\ 
- \Delta_4(3;\mu;M_1)_i \\ \end{array}\right) = 
\left(\begin{array}{cccc} 
- \Delta_1(2;\mu_1)_{i_1} & \Delta_1(2;\mu_2)_{i_2}  
& - \Delta_1(2;\mu_3)_{i_3} & 0 \\  
\Delta_2(2;\mu_1)_{i_1} & - \Delta_2(2;\mu_2)_{i_2} & 0   
& \Delta_1(2;\mu_2)_{i_4} \\  
- \Delta_3(2;\mu_1)_{i_1} & 0 & \Delta_2(2;\mu_3)_{i_3}  
& - \Delta_2(2;\mu_2)_{i_4} \\ 
0 & \Delta_3(2;\mu_2)_{i_2} & - \Delta_3(2;\mu_3)_{i_3}  
& \Delta_3(2;\mu_2)_{i_4} \\ \end{array}\right) \cdot 
\left(\begin{array}{c}  x_1 \\ x_2 \\ x_3 \\ x_4 \\ 
\end{array}\right)\ .$$

En rajoutant~: 
$$\begin{array}{l} R(3;m_1,m_2,m_3;M)_i = - x_1\cdot \biggl\{\Delta_1(2;\mu_1)_{i_1} -  
\Delta_2(2;\mu_1)_{i_1} + \Delta_3(2;\mu_1)_{i_1}\biggr\} \\[1em]\ \ \ 
+ x_2\cdot \biggl\{\Delta_1(2;\mu_2)_{i_2} -  
\Delta_2(2;\mu_2)_{i_2} + \Delta_3(2;\mu_2)_{i_2}\biggr\} 
- x_3\cdot \biggl\{\Delta_1(2;\mu_3)_{i_3} -  
\Delta_2(2;\mu_3)_{i_3} + \Delta_3(2;\mu_3)_{i_3}\biggr\} \\[1em]\ \ \
 + x_4 \cdot \biggl\{\Delta_1(2;\mu_2)_{i_4} - 
\Delta_2(2;\mu_2)_{i_4} + \Delta_3(2;\mu_2)_{i_4}\biggr\} = 0\ .\end{array}$$

Le th\'eor\`eme ci-dessous g\'en\'eralise ces exemples. 

\medskip
{\bf 2.7.} {\it Th\'eor\`eme.} On a 
$$ R(n;m_1,\ldots,m_n;M)_i = 0 $$
pour tous $n, m_1, \ldots, m_n, M$ et $i$. 

{\it D\'emonstration}~: elle se fait par r\'ecurrence sur $n$. Le cas $n=2$ est 
l'exemple 2.5.  

Le passage de $n-1$ \`a $n$ suit l'exemple 2.6. 
  
Posons pour abr\'eger 
$$\mu = (m_1, \ldots, m_n)\ .$$
\`A partir de cela, on introduit $n+1$ vecteurs 
$\mu_j \in \BZ^{n-1}$~: 
$$\mu_j := (m_1, \ldots , \hat m_j, \ldots, m_n), $$
$j = 1, \ldots, n$, et 
$$\mu_{n+1} := (m_1, \ldots, \hat m_{n-1}, m_n) = \mu_{n-1}\ .$$
On d\'efinit le vecteur 
$$(i_1,i_2, \ldots, i_{n+1}) := (\underbrace{i, i, \ldots, i}_{n-1\ {\mathrm{fois}}}, 
i + m_{n-1} - m_n, m_{n-1} + m_n - 1) \in \BZ^{n+1}\ .$$
En d\'eveloppant les d\'eterminants $\Delta_j(n;\mu,M_{n+2-j})_i$, $2 \leq j \leq n+1$ 
suivant la deuxi\`eme colonne et le d\'eterminant $\Delta_1(n;\mu,M_{n+1})_i$ suivant la 
deuxi\`eme ligne, on obtient~: 
$$R(n;\mu;M)_i = \sum_{j=1}^{n+1}\ (-1)^j x_j R(n-1;\mu_j;M_{1j})_{i_j}\ .$$
Ici $M_{1j}$ est la matrice obtenue en enlevant la premi\`ere ligne et 
la $j$-i\`eme colonne de la matrice $M$. 

Notre assertion en d\'ecoule imm\'ediatement par r\'ecurrence sur $n$. 

\medskip
{\bf 2.8.} On aura besoin d'un cas particulier de ces relations. Prenons 
$$\mu = (m_1,m_2,\ldots,m_n) = (2,3,\ldots, n+1)\ .$$
Pour la matrice $M$, prenons
$$M = \left(\begin{array}{cccccc} b(1)_3 & b(2)_4 & b(2)_5 & ... & b(2)_{n+2} & b(2)_{i+n+1} \\ 
b(1)_4 & b(2)_5 & b(3)_5 & ... & b(3)_{n+3} & b(3)_{i+n+2} \\ 
. & . & . & . & . & . \\ 
b(1)_{n} & b(2)_{n+1} & b(3)_{n+2} & \ldots & b(n-1)_{2n-1} & b(n-1)_{i+2n-2} \\ 
\end{array}\right)\ .$$
Alors le premier d\'eterminant 
$$\Delta_1(n;\mu;M_{n+1})_{i+2n} = c(n+1)_i\ .$$
On pose par d\'efinition~: 
$$c(n+1)'_i := \Delta_2(n;\mu;M_{n})_{i+2n}\ ,$$
$$c(n+1)''_i := \Delta_3(n;\mu;M_{n-1})_{i+2n}\ .$$
Par contre, si $j \geq 4$ on voit que dans le d\'eterminant 
$\Delta_j(n;\mu;M_{n+2-j})_{i+2n}$ la derni\`ere colonne est \'egale \`a la 
$(n-j+3)$-i\`eme colonne, d'o\`u 
$$ \Delta_4(n;\mu;M_{n-2})_{i+2n} = \Delta_5(n;\mu;M_{n-3})_{i+2n} = 
\ldots = \Delta_{n+1}(n;\mu;M_{1})_{i+2n} = 0\ .$$
Donc 2.7 entra\^ine 

\medskip
{\bf 2.9.} {\it Corollaire.} Pour tous $n \geq 3$ 
$$ c(n)_i - c(n)'_i + c(n)''_i = 0\ .$$    
       
\bigskip\bigskip

\centerline{\bf \S\  3. D\'ebut de la d\'emonstration du th\'eor\`eme 1.5}\label{Debut} 

\bigskip

{\bf 3.1.} On a  
$$ f(x) = a_nx^n + a_{n-1}x^{n-1} + \ldots + a_0\ .$$
La d\'eriv\'ee~: 
$$\begin{array}{l} f_1(x) = f'(x) = na_nx^{n-1} + (n-1)a_{n-1}x^{n-2} + \ldots + a_1 
= n a_n\biggl\{x^{n-1} +  \frac{(n-1)a_{n-1}}{na_n}x^{n-2} + \ldots + 
\frac{a_1}{na_n}\biggr\} \\[1em]
\ \ \ \ \ \ \ \  = \gamma_1(c(1)_0 x^{n-1} + c(1)_1 x^{n-2} + \ldots + c(1)_{n-1})\ .\end{array}$$ 

\medskip
{\bf 3.2.} Le quotient de la division euclidienne de deux polyn\^omes $f(x)$ et
$g(x)=a'_{n-1}x^{n-1} +a'_{n-2}x^{n-2}+\dots$ est \'egal \`a 
$$ \textstyle\frac{a_n}{a'_{n-1}}\, x+{a'_{n-1}a_{n-1}-a'_{n-2}a_n\over (a'_{n-1})^2}\ .$$

On fait la division euclidienne~: 
%$$ f - (x/n)f' = (a_n/n)x^{n-1} + (2a_{n-1}/n)x^{n-2} + (3a_{n-2}/n)x^{n-3} + 
%\ldots, $$
%ensuite
$$\begin{array}{l} f - (x/n + a_{n-1}/n^2a_n)f' 
 = \frac{2n a_na_{n-2} - (n-1)a_{n-1}^2}{n^2a_n}x^{n-2} +  
\frac{3n a_na_{n-3} - (n-2)a_{n-1}a_{n-2}}{n^2a_n}x^{n-3} + \ldots\\[1em]
\ \ \ \ \  = \frac{1}{n^2a_n}\cdot (b(1)_2x^{n-2} + b(1)_3x^{n-3} + \ldots + b(1)_n) %\\[1em]
 = - \gamma_2\cdot (c(2)_0 x^{n-2} + c(2)_1 x^{n-3} + \ldots c(2)_{n-2})\ .\end{array}$$ 
Donc 
$$ f_2(x) = \gamma_2\cdot \sum_{i=0}^{n-2}\ c(2)_i x^{n-2-i}\ .$$
Cela d\'emontre l'assertion 1.5 pour $j= 1, 2$, et l'on proc\`ede par
r\'ecurrence par $j$. 

\medskip
{\bf 3.3.} On suppose que l'on a d\'ej\`a trouv\'e~: 
$$ f_{j-1}(x) = \gamma_{j-1}\cdot [c(j-1)x^{n-j+1} + c(j-1)_1x^{n-j} + \ldots 
+ c(j-1)_ix^{n-j+1-i} + \ldots ] $$
et
$$ f_j(x) = \gamma_j\cdot [c(j)x^{n-j} + c(j)_1x^{n-j-1} + \ldots 
+ c(j)_ix^{n-j-i} + \ldots ] \ . $$
On fait la division euclidienne~: 
%$$ f_{j-1}(x) - \frac{\gamma_{j-1}c(j-1)}{\gamma_jc(j)} x\cdot f_j(x) = 
%\gamma_{j-1}\biggl[c(j-1)_1 - \frac{c(j-1)}{c(j)}\cdot c(j)_1\biggr]
%\cdot x^{n-j} + $$
%$$ + \ldots + 
%\gamma_{j-1}\biggl[c(j-1)_i - \frac{c(j-1)}{c(j)}\cdot c(j)_i\biggr]
%\cdot x^{n-j+1-i} 
%+ \ldots $$
%Ensuite,
$$ f_{j-1}(x) - \Bigl(\frac{\gamma_{j-1}c(j-1)}{\gamma_jc(j)} x +
 \gamma_{j-1}\biggl[c(j-1)_1 - \frac{c(j-1)}{c(j)}\cdot c(j)_1\biggr]\cdot 
\frac{1}{\gamma_j c(j)}\Bigr) f_j(x) = $$
$$ = \sum_{i=2}^{n-j+1}\ \gamma_{j-1}\biggl\{c(j-1)_i - \frac{c(j-1)}{c(j)}
\cdot c(j)_i 
- \biggl[c(j-1)_1 - \frac{c(j-1)}{c(j)}\cdot c(j)_1\biggr]\cdot 
\frac{c(j)_{i-1}}{c(j)}\biggr\}\cdot x^{n-j+1-i} = $$
$$ = \frac{\gamma_{j-1}}{c(j)^2}\sum_{i=2}^{n-j+1}\ 
\biggl\{c(j-1)_ic(j)^2 - c(j-1)c(j)c(j)_i 
- c(j-1)_1c(j)_{i-1}c(j) + c(j-1)c(j)_1c(j)_{i-1}\biggr\}\cdot x^{n-j+1-i}\ .$$
On pose~: 
$$ Q(j)_i := c(j-1)_ic(j)^2 - c(j-1)c(j)c(j)_i 
- c(j-1)_1c(j)_{i-1}c(j) + c(j-1)c(j)_1c(j)_{i-1}
\eqno{(3.3.1)} $$
Alors on a~: 
$$ f_{j+1}(x) = - \frac{\gamma_{j-1}}{c(j)^2}\sum_{i=2}^{n-j+1}\ 
Q(j)_i x^{n-j+1-i} $$
Il faut montrer que 
$$ f_{j+1}(x) = \gamma_{j+1} \sum_{i=0}^{n-j-1}\ c(j+1)_ix^{n-j-1-i} = 
 \gamma_{j+1} \sum_{i=0}^{n-j+1}\ c(j+1)_{i-2}x^{n-j+1-i} $$
o\`u
$$ \gamma_{j+1} = \gamma_{j-1}\cdot \frac{c(j-1)^2}{c(j)^2}\ .$$
Donc notre th\'eor\`eme est \'equivalent \`a l'identit\'e suivante~: 
$$ Q(j)_i = - c(j)^2c(j+1)_i\ .\eqno{(3.3.2)}$$
 
\bigskip\bigskip  

\centerline{\bf \S\  4. Formule $(A)$}\label{FormuleA} 

\bigskip 

{\bf 4.1.} Revenons \`a notre alg\`ebre $\fB$.  

On consid\`ere la matrice $n\times n$
$$ C(n+1)_{i-2} = \left(\begin{array} {ccccc}
b(1)_2 & b(1)_3 & \ldots & b(1)_n & b(1)_{n+1} \\ 
b(1)_3 & b(2)_4 & \ldots & b(2)_{n+1} & b(2)_{n+2} \\ 
. & . & \ldots & . & . \\ 
b(1)_n & b(2)_{n+1} & \ldots & b(n-1)_{2n-2} & b(n-1)_{2n-1} \\ 
b(1)_{n+i-1} & b(2)_{n+i} & \ldots & b(n-1)_{2n+i-3} & b(n)_{2n+i-2} \\
\end{array}\right)\ .$$
Donc $c(n+1)_{i-2} = \det C(n+1)_{i-2}$. 

Si l'on d\'esigne par $C(n+1)_{i-2;\hat p, \hat q}$ la matrice $C(n+1)_{i-2}$ 
avec la $p$-i\`eme ligne et la $q$-i\`eme colonne enlev\'ee, on aura~: 
$$ c(n) = \det C(n+1)_{i-2;\hat n, \hat n}\ ,$$
$$  c(n)_{i-1} = \det C(n+1)_{i-2;\hat {n-1}, \hat n}\ .$$
En plus, on a~: 
$$ c(n)''_i = \det C(n+1)_{i-2;\hat {n-2}, \hat n}  $$
o\`u $c(n)''_i$ a \'et\'e introduit dans 2.9. 

\medskip
{\bf 4.2.} {\it Th\'eor\`eme.} Pour tous $n, i \in \BZ,\ n\geq 3,$ on a 
la relation suivante dans $\fB$
$$ c(n-1)_ic(n)^2 - c(n-1)c(n)c(n)_i - c(n-1)_1c(n)_{i-1}c(n) 
+ c(n-1)c(n)_1c(n)_{i-1} $$
$$ = - c(n-1)^2 c(n+1)_{i-2}\eqno{(F)} $$

On a vu que notre th\'eor\`eme principal 1.5 est une cons\'equence de $(F)$~: 
en effet $(F)$ co\"incide avec la formule (3.3.2) (avec $j$ remplac\'e par 
$n$).  

\`A son tour, $(F)$ est une cons\'equence imm\'ediate de deux formules~: 
$$ c(n-1)_ic(n) - c(n-1)_1c(n)_{i-1} = - c(n-1)c(n)''_i
\eqno{(A)} $$
ou bien 
$$ c(n-1)_ic(n) - c(n-1)_1c(n)_{i-1} + c(n-1)c(n)''_i = 0\eqno{(A')} $$
et 
$$ \bigl\{c(n)_i + c(n)''_i\bigr\}\cdot c(n) - 
c(n)_1 c(n)_{i-1} = c(n-1) c(n+1)_{i-2}\ .\eqno{(B)}$$

La d\'emonstration de (B) utilise les relations quadratiques entre les lettres 
$b(i)_j$. Par contre, (A) est "\'el\'ementaire", en ce sens que cette 
identit\'e n'utilise pas de relations entre les lettres $b(i)_j$.  

Pour d\'emontrer (A), on applique le lemme suivant (une variante des
relations de Pl\"ucker)~:  
  
\medskip
{\bf 4.3.} {\it Lemme.} $(A_n)$\  Consid\'erons $n$ vecteurs de dimension $n-1$, 
$w_i = (w_{i1},\ldots,w_{i,n-1}),\ i= 1, \ldots , n$. \`A partir d'eux, 
on d\'efinit $n$ vecteurs de dimension $n-2$~:  
$v_i = (w_{i1},\ldots, w_{i,n-2})$. On pose~: 
$$ W_i = \det(w_1,\ldots,\hat w_i,\ldots, w_n)^t\ ,$$
$$ V_{ij}:= \det(v_1,\ldots, \hat v_i,\ldots, \hat v_j, \ldots, v_{n})^t\ .$$
Alors 
$$ V_{n-2,n-1}\cdot W_{n} - V_{n-2,n}\cdot W_{n-1} + 
V_{n-1,n}\cdot W_{n-2} = 0\ .$$ 
  
$(B_n)$\  Consid\'erons $n$ vecteurs de dimension $n-2$, 
$v_i = (v_{i1},\ldots,v_{i,n-2}),\ i= 1, \ldots , n$. Consid\'erons les 
mineurs 
$$ V_{ij}:= \det(v_1,\ldots, \hat v_i,\ldots, \hat v_j, \ldots, v_{n})^t\ .$$
Alors pour chaque $i < n-2$,  
$$ V_{n-2,n-1}\cdot V_{i, n} - V_{n-2,n}\cdot V_{i, n-1} + 
V_{n-1,n}\cdot V_{i, n-2} = 0\ .$$

En effet, en d\'eveloppant $W_i$ par rapport \`a la derni\`ere colonne, 
on obtient~: $(B_n) \Rightarrow (A_n)$. 

Par contre, pour v\'erifier $(B_n)$, consid\'erons la matrice $(n-1)\times 
(n-2)$, $W^\sim = V_i$.  Alors on aura $V_{ij} = W^\sim_j,\ j = n, n-1, n-2$. 
D'un autre c\^ot\'e, en d\'eveloppant les mineurs dans $(B_n)$~: 
$V_{pq},\ n-2\leq p < q \leq n$ par rapport \`a la $i$-i\`eme ligne, on 
obtient les mineurs $V^\sim_{pq}$, o\`u $V^\sim$ est obtenue de $W^\sim$ 
en enlevant la derni\`ere colonne. On v\'erifie que $(B_n)$ se r\'eduit 
\`a $(A_{n-1})$ correspondant \`a $W^\sim$. 

Il s'ensuit que $(A_{n-1}) \Rightarrow (B_n)$ et on conclut par r\'ecurrence. 

\medskip
{\bf 4.4.} Le lemme \'etant v\'erifi\'e, l'assertion 4.2 $(A)$ est 4.3 $(A_n)$ 
pour la matrice $W$ \'egale \`a $c(n+1)_{i-2}$ avec la derni\`ere colonne 
enlev\'ee.      

\bigskip\bigskip 

\centerline{\bf \S\  5. Formule $(B)$}\label{FormuleB} 

\bigskip 

{\bf 5.1.} Maintenant on s'occupe de la formule 
$$ P:= \bigl\{c(n)_i + c(n)''_i\bigr\}\cdot c(n) - 
c(n)_1 c(n)_{i-1} = c(n-1) c(n+1)_{i-2}\ .
\eqno{(B)} $$
On introduit $n$ vecteurs de dimension $n-1$, 
$w_1, \ldots, w_n$ qui sont les lignes de la matrice $c(n+1)_{i-2}$ 
sans la derni\`ere colonne~:   
$$ W = \left(\begin{array}{ccccc} 
b(1)_2 & b(1)_3 & \ldots & b(1)_{n-1} & b(1)_n  \\ 
b(1)_3 & b(2)_4 & \ldots & b(2)_{n} & b(2)_{n+1}  \\ 
. & . & \ldots & . & . \\
b(1)_{n-1} & b(2)_{n} & \ldots & b(n-2)_{2n-4} & b(n-2)_{2n-3}  \\ 
b(1)_n & b(2)_{n+1} & \ldots & b(n-2)_{2n-3} & b(n-1)_{2n-2}  \\ 
b(1)_{n+i-1} & b(2)_{n+i} & \ldots & b(n-2)_{2n+i-4} & b(n-1)_{2n+i-3}  \\
\end{array}\right) $$
et $n$ mineurs
$$ W_i = \det(w_1,\ldots,\hat w_i,\ldots, w_n)^t,\ 
i = 1, \ldots, n\ .$$
Par exemple, $W_n = c(n),\ W_{n-1} = c(n)_{i-1},\ W_{n-2} = c(n)''_i$. 
Donc,  
$$\begin{array}{l} c(n+1)_{i-2} = b(n)_{2n+i-2}W_n - b(n-1)_{2n-1}W_{n-1} \\[1em]
\ \ \ \ \ \ \ \ \ \ \ \ \ \ \  + b(n-2)_{2n-2}W_{n-2} - \ldots + (-1)^{n-1}b(1)_{n+1}W_1\\[1em] 
\ \ \ \ \ \ \ \ \ \ \ \ \ \ \  = b(n)_{2n+i-2}W_n - b(n-1)_{2n-1}W_{n-1} + R\end{array}\eqno{(5.1.1)} $$
o\`u
$$ R = b(n-2)_{2n-2}W_{n-2} - b(n-3)_{2n-3}W_{n-3} + 
\ldots + (-1)^{n-1}b(1)_{n+1}W_1\ .\eqno{(5.1.2)}$$

\medskip
{\bf 5.2.} On a $n-1$ relations lin\'eaires entre les $W_i$~: la $i$-i\`eme est 
obtenue en ajoutant \`a $W$ sa $i$-i\`eme colonne et en d\'eveloppant 
le d\'eterminant $=0$ par rapport \`a la derni\`ere colonne. 

Explicitement~:
$$ b(n-1)_{2n+i-3}W_n - b(n-1)_{2n-2}W_{n-1} + 
b(n-2)_{2n-3}W_{n-2} - \ldots + (-1)^{n-1}b(1)_{n}W_1 = 0\ ,$$
$$ b(n-2)_{2n+i-4}W_n - b(n-2)_{2n-3}W_{n-1} + 
b(n-2)_{2n-4}W_{n-2} - \ldots + (-1)^{n-1}b(1)_{n-1}W_1 = 0\ ,$$
$$ .\ .\ .\ $$
$$ b(2)_{n+i}W_n - b(2)_{n+1}W_{n-1} + 
b(2)_{n}W_{n-2} - \ldots + (-1)^{n-2}b(2)_4W_2 + (-1)^{n-1}b(1)_{3}W_1 = 0\ ,$$   
$$ b(1)_{n+i-1}W_n - b(1)_{n}W_{n-1} + 
b(1)_{n-1}W_{n-2} - \ldots + (-1)^{n-2}b(1)_3W_2 + (-1)^{n-1}b(1)_{2}W_1 = 0\ .$$

\medskip
{\bf 5.3.} D'autre part, rappelons la matrice $c(n)_1$~:
$$ c(n)_1 = \det\left(\begin{array}{ccccc} 
b(1)_2 & b(1)_3 & \ldots & b(1)_{n-1} & b(1)_n  \\ 
b(1)_3 & b(2)_4 & \ldots & b(2)_{n} & b(2)_{n+1}  \\ 
. & . & \ldots & . & . \\
b(1)_{n-1} & b(2)_{n} & \ldots & b(n-2)_{2n-4} & b(n-2)_{2n-3}  \\  
b(1)_{n+1} & b(2)_{n+2} & \ldots & b(n-2)_{2n-2} & b(n-1)_{2n-1}  \\ 
\end{array}\right)\ .$$
On d\'eveloppe cette quantit\'e par rapport \`a la derni\`ere colonne~: 
$$ c(n)_1  = b(n-1)_{2n-1}c(n-1) - b(n-2)_{2n-3}M_{n-2} + 
\ldots + (-1)^{n-1}b(2)_{n+1}M_2 + (-1)^{n}b(1)_{n}M_1\ .$$
Apr\`es la multiplication par $- c(n)_{i-1} = - W_{n-1}$ on obtient~: 
$$ - c(n)_1c(n)_{i-1} = - b(n-1)_{2n-1}c(n-1)W_{n-1}\ (*) + R' $$
o\`u
$$ R' = b(n-2)_{2n-3}W_{n-1}M_{n-2} - b(n-3)_{2n-4}W_{n-1}M_{n-3} + \ldots $$
$$ + (-1)^{n}b(2)_{n+1}W_{n-1}M_2 + (-1)^{n-1}b(1)_{n}W_{n-1}M_1\ .$$

\medskip
{\bf 5.4.} 
Maintenant rempla\c{c}ons dans $R'$ les termes $(-1)^ib(n-i)_{2n-i-1}W_{n-1}$ 
en utilisant les relations 5.2~:   
$$ b(n-2)_{2n-3}W_{n-1} = b(n-2)_{2n+i-4}W_n  + 
b(n-2)_{2n-4}W_{n-2} - \ldots + (-1)^{n-1}b(1)_{n-1}W_1\ ,$$
$$ .\ .\ .\ $$
$$ b(2)_{n+1}W_{n-1} =  b(2)_{n+i}W_n  + 
b(2)_{n}W_{n-2} - \ldots + (-1)^{n-2}b(2)_4W_2 + (-1)^{n-1}b(1)_{3}W_1\ ,$$   
$$ b(1)_{n}W_{n-1} =  b(1)_{n+i-1}W_n  + 
b(1)_{n-1}W_{n-2} - \ldots + (-1)^{n-2}b(1)_3W_2 + (-1)^{n-1}b(1)_{2}W_1\ .$$
Alors on obtient~:
$$ - c(n)_1c(n)_{i-1} = - b(n-1)_{2n-1}c(n-1)W_{n-1}\ (*) $$
$$ + \bigl\{b(n-2)_{2n+i-4}M_{n-2} - \ldots + (-1)^nb(2)_{n+i}M_2 
+ (-1)^{n+1}b(1)_{n+i-1}M_1\bigr\}\cdot c(n) + R'', $$
o\`u~: 
$$ R'' = \biggl\{b(n-2)_{2n-4}W_{n-2} - \ldots + (-1)^{n-1}b(1)_{n-1}W_1 \biggr\} 
\cdot M_{n-2} - \ldots $$
$$ + (-1)^n\cdot \biggl\{
b(2)_{n}W_{n-2} - \ldots + (-1)^{n-2}b(2)_4W_2 + 
(-1)^{n-1}b(1)_{3}W_1 \biggr\}\cdot M_2 $$
$$ + (-1)^{n-1}\cdot\biggl\{
b(1)_{n-1}W_{n-2} - \ldots + (-1)^{n-2}b(1)_3W_2 + (-1)^{n-1}b(1)_{2}W_1\biggr\}
\cdot M_1\ .$$

\medskip
{\bf 5.5.} {\it Lemme.} $R'' = c(n-1)R$. 

{\it D\'emonstration.} 
On introduit les vecteurs de dimension $n - 2$~: 
$$ \CW = \bigl((-1)^{n+1} W_1, (-1)^{n+2} W_2, \ldots, W_{n-2}\bigr)\ ,$$  
$$ \CM = \bigl((-1)^{n+1} M_1, (-1)^{n+2} M_2, \ldots, M_{n-2}\bigr) $$
et
$$ b = \bigl(b(1)_{n+1}, b(1)_{n+2}, \ldots, b(1)_{2n-2}\bigr)\ .$$ 
Alors la d\'efinition de $R''$ se r\'ecrit~: 
$$ R'' = \CM\cdot C(n-1) \cdot \CW^t
\eqno{(5.5.1)} $$
(o\`u $c(n-1) = \det C(n-1)$, la matrice $C(n-2)$ \'etant sym\'etrique)~; 
de plus,
$$ R = b\cdot \CW^t\ .$$ 
Maintenant d\'eveloppons les quantit\'es 
$M_i$ par rapport \`a la derni\`ere ligne~: 
$$ M_i = b(1)_{2n-2} M_{i,n-2} - b(1)_{2n-3} M_{i,n-3} + \ldots 
+ (-1)^{n+2}b(1)_{n+2}M_{i2} + (-1)^{n+1}b(1)_{n+1}M_{i1} $$
$$ = \sum_{j=1}^{n-2}\ (-1)^{n+j}b(1)_{n+j} M_{ij},\ i=1, \ldots, n-2\ . $$
On remarque que les quantit\'es  $M_{ij}$ sont les mineurs de la matrice 
$(n-2)\times (n-2)$ $C(n-1)$. Il vient~: 
$$ \CM = b \cdot \hat C(n-1) $$
o\`u 
$$ \hat C(n-1) = \bigl( (-1)^{i+j} M_{ij}\bigr), $$
donc $\hat C(n-1)\cdot C(n-1) = c(n-1)$. En substituant dans (5.5.1)~: 
$$ R'' = b\cdot \hat C(n-1)\cdot C(n-1) \cdot \CW^t = c(n-1)\cdot b\cdot \CW^t 
= c(n-1)R, $$   
cqfd. 

\medskip
{\bf 5.6.} Il s'ensuit que pour v\'erifier l'identit\'e (B) il reste \`a 
d\'emontrer que 
$$\begin{array}{l} c(n)_i + c(n)''_i + b(n-2)_{2n+i-4}M_{n-2} - \ldots \\[1em]
\ \ \ \ \ \ \ \ \ \  + (-1)^nb(2)_{n+i}M_2 
+ (-1)^{n+1}b(1)_{n+i-1}M_1 = b(n)_{2n+i-2}c(n-1)\ .\end{array} 
\eqno{(5.6.1)} $$
Par contre, la quantit\'e 
$$ b(n)_{2n+i-2}c(n-1) - 
b(n-2)_{2n+i-4}M_{n-2} + \ldots + 
(-1)^{n+1}b(2)_{n+i}M_2 
+ (-1)^{n}b(1)_{n+i-1}M_1 $$ 
n'est autre que le d\'eveloppement de $c(n)'_i$ suivant la derni\`ere colonne, donc 
(5.6.1) est \'equivalent \`a 
$$ c(n)_i + c(n)''_i = c(n)'_i
\eqno{(5.6.2)} $$
qui a \'et\'e d\'ej\`a prouv\'ee, cf. Corollaire 2.9. 

Ceci ach\`eve la d\'emonstration du th\'eor\`eme 4.2, et donc du 1.5. 

\medskip
{\bf 5.7.} {\it D\'emonstration du th\'eor\`eme} 1.9. En fait, nous l'avons 
d\'ej\`a montr\'e~: la d\'emonstration de la r\'ecurrence principale 
(3.3.2) n'utilise que les relations dans l'alg\`ebre $\fB$. 

Ces relations 
sont v\'erifi\'ees si l'on d\'efinit les variables $b(i)_j$ \`a partir de 
coefficients de polyn\^omes $f_1(x)$ et $f_2(x)$ comme dans 1.8, 
d'o\`u l'assertion.     

\bigskip\bigskip\bigskip

\newpage 
\centerline{DEUXI\`EME PARTIE.} 

\bigskip\bigskip

\centerline{POLYN\^OMES D'EULER ET D\'ETERMINANT DE CAUCHY}

\bigskip\bigskip

\centerline{\bf \S\  1. Nombres $\beta(j)_i$}\label{Nombresbeta}

\bigskip 

{\bf 1.1.} 
Rappelons que pour un polyn\^ome 
$$ f(x) = a_n x^n + a_{n-1}x^{n-1} + \ldots + a_0 $$
les nombres $b(j)_i$ sont d\'efinis par 
$$ b(j)_i = n\sum_{p=0}^{j-1}\ (i-2p)a_{n-p}a_{n-i+p} 
- j(n-i+j)a_{n-j}a_{n+j-i}\ .$$
On introduit les quantit\'es~: 
$$ q_i := \frac{a_{i-1}}{a_i}\ ,$$
$$ r_i := \frac{q_{i-1}}{q_i} = \frac{a_ia_{i-2}}{a_{i-1}^2}\ ,  $$
puis
$$ \beta(j)_i := \frac{b(j)_i}{(n-i+j)a_{n-j}a_{n+j-i}} = \sum_{p=0}^{j-1}\ \frac{n(i-2p)}{n-i+j}\cdot 
\frac{a_{n-p}a_{n+p-i}}{a_{n-j}a_{n+j-i}} - j\ .$$

\medskip
{\bf 1.2.} Par exemple~: 
$$ \beta(1)_2 = \frac{2n}{n-1}\cdot\frac{a_{n}a_{n-2}}{a_{n-1}^2} - 1 = 
\frac{2n}{n-1}\cdot r_n - 1\ ,$$
$$ \beta(1)_i = \frac{ni}{n-i+1}\cdot 
\frac{a_{n}a_{n-i}}{a_{n-1}a_{n-i+1}} - 1\ .$$
On remarque que
$$\frac{a_{n}a_{n-i}}{a_{n-1}a_{n-i+1}} = \frac{q_{n-i+1}}{q_n} = 
r_{n-i+2}r_{n-i+1}\ldots r_n\ .$$
On d\'efinit les quantit\'es
$$\psi(i,j):= \prod_{p=i}^j\ r_p$$
(donc $\psi(i,j) = 1$ si $i > j$). Il s'ensuit~: 
$$\beta(1)_i = \frac{ni}{n-i+1}\cdot \psi(n-i+2,n) - 1\ .$$

\medskip
{\bf 1.3.} De m\^eme~: 
$$ \frac{a_{n}a_{n-i}}{a_{n-2}a_{n-i+2}} = \frac{a_{n}a_{n-i}}{a_{n-1}a_{n-i+1}} \cdot 
\frac{a_{n-1}a_{n-i+1}}{a_{n-2}a_{n-i+2}} = \psi(n-i+2,n)\psi(n-i+3,n-1)\ .$$
Par exemple~: 
$$ \frac{a_{n}a_{n-4}}{a_{n-2}^2} = \psi(n-2,n)\psi(n-1,n-1) = r_{n-2}r_{n-1}^2r_n\ .$$
Il en d\'ecoule~:  
$$ \beta(2)_4 = \frac{4n}{n-2}\frac{a_na_{n-4}}{a_{n-2}^2} 
+ \frac{2n}{n-2}\frac{a_{n-1}a_{n-3}}{a_{n-2}^2} - 2  = \frac{4n}{n-2}r_{n-2}r_{n-1}^2r_n + \frac{2n}{n-2}r_{n-1} - 2\ ,$$

$$\beta(2)_i = \frac{ni}{n-i+2}\frac{a_n a_{n-i}}{a_{n-2}a_{n-i+2}} 
+ \frac{n(i-2)}{n-i+2}\frac{a_{n-1} a_{n-i+1}}{a_{n-2}a_{n-i+2}} - 2 $$
$$ = \frac{ni}{n-i+2}\psi(n-i+2,n)\psi(n-i+3,n-1)+\frac{n(i-2)}{n-i+2}\psi(n-i+3,n-1) - 2\ .$$

\medskip
{\bf 1.4.} Un autre exemple~: 
$$\frac{a_n a_{n-6}}{a_{n-3}^2} = \psi(n-4,n)\psi(n-3,n-1)\psi(n-2,n-2)
= r_{n-4}r_{n-3}^2r_{n-2}^3r_{n-1}^2r_{n}\ .$$

\medskip
{\bf 1.5.} En g\'en\'eral on pose~: 
$$ \phi(n,j,i):= \frac{a_na_{n-i}}{a_{n-j}a_{n-i+j}} = \prod_{q=0}^{j-1}\ \psi(n-i+j+q,n-q)$$
et l'on aura~:
$$ \beta(j)_i = \sum_{p=0}^{j-1}\ \frac{n(i-2p)}{n-i+j}\cdot \phi(n-p,j-p,i-p) - j\ .$$

\medskip
{\bf 1.6.} Passons maintenant aux d\'eterminants $c(n)$. On commence 
par un exemple~: 
$$ c(4) = \det\left(\begin{array}{ccc}
b(1)_2 & b(1)_3 & b(1)_4 \\ 
b(1)_3 & b(2)_4 & b(2)_5 \\
b(1)_4 & b(2)_5 & b(3)_6 \\ 
\end{array}\right) $$

$$ = \det\left(\begin{array}{ccc} 
(n-1)a_{n-1}^2\beta(1)_2 & (n-2)a_{n-1}a_{n-2}\beta(1)_3 
& (n-3)a_{n-1}a_{n-3}\beta(1)_4 \\ 
(n-2)a_{n-1}a_{n-2}\beta(1)_3 & (n-2)a_{n-2}^2\beta(2)_4 
& (n-3)a_{n-2}a_{n-3}\beta(2)_5 \\
(n-3)a_{n-1}a_{n-3}\beta(1)_4 & (n-3)a_{n-2}a_{n-3}\beta(2)_5 
& (n-3)a_{n-3}^2\beta(3)_6 \\ 
\end{array}\right) $$

$$ = (a_{n-1}a_{n-2}a_{n-3})^2\cdot 
\det\left(\begin{array}{ccc} 
(n-1)\beta(1)_2 & (n-2)\beta(1)_3 
& (n-3)\beta(1)_4 \\ 
(n-2)\beta(1)_3 & (n-2)\beta(2)_4 
& (n-3)\beta(2)_5 \\
(n-3)\beta(1)_4 & (n-3)\beta(2)_5 
& (n-3)\beta(3)_6 \\ 
\end{array}\right)\ .$$

\medskip
{\bf 1.7.} En g\'en\'eral
$$ c(m+1) = \biggl(\prod_{i=1}^{m}\ a_{n-i}\biggr)^2 \times  $$

$$ \times \det\left(\begin{array}{cccc} 
(n-1)\beta(1)_2 & (n-2)\beta(1)_3  & \ldots 
& (n-m)\beta(1)_{m+1} \\ 
(n-2)\beta(1)_3 & (n-2)\beta(2)_4  &  \ldots 
& (n-m)\beta(2)_{m+2}\\
 . & . & \ldots & . \\ 
(n-m)\beta(1)_{m+1} & (n-m)\beta(2)_{m+2} & 
\ldots & (n-m)\beta(m)_{2m} \\  
\end{array}\right)\ .$$   

\bigskip\bigskip\bigskip

\centerline{\bf \S\  2. Polyn\^omes d'Euler et fonction hyperg\'eom\'etrique}\label{PolynomesEuler} 

\bigskip

{\bf 2.1.}  Suivant [Euler], on d\'efinit 
les polyn\^omes 
$$ E_n(x) = \frac{1}{2}\{(1 + ix/{2n})^{2n} + (1 - ix/{2n})^{2n}\}\ .
\eqno{(2.1.1)} $$
Donc, $E_{n}(x)$ est un polyn\^ome de degr\'e $2n$, avec le terme 
constant $1$, ne contenant que des puissances paires de $x$. Plus pr\'ecis\'ement, 
$$ E_{n}(x) = \sum_{k=0}^{n}\ (-1)^k\binom{2n}{2k}\frac{x^{2k}}{(2n)^{2k}}\ .
\eqno{(2.1.2)}$$ 
Par exemple~: 
$$ E_1(x) = 1 - \frac{1}{4}x^2\ ,$$
$$ E_2(x) = 1 - \frac{3}{8}x^2 + \frac{1}{256}x^4\ ,$$
$$ E_3(x) = 1 - \frac{5}{12}x^2 + \frac{5}{432}x^4 - \frac{1}{46656}x^4\ ,$$
$$ E_4(x) =  1 - \frac{7}{16}x^2 + \frac{35}{2048}x^4 - \frac{7}{65536}x^6 
+ \frac{1}{16777216}x^8\ .$$

\medskip
{\bf 2.2.} Rappelons que la fonction hyperg\'eom\'etrique de Gauss 
est d\'efinie par
$$ F(\alpha,\beta,\gamma,x) = 1 + \frac{\alpha\beta}{1\cdot\gamma}x + 
\frac{\alpha(\alpha +1)\beta(\beta + 1)}{1\cdot 2\cdot\gamma(\gamma + 1)}x^2 + 
\frac{\alpha(\alpha + 1)(\alpha + 2)\beta(\beta + 1)(\beta + 2)} 
{1 \cdot 2 \cdot 3 \cdot \gamma(\gamma + 1)(\gamma + 2)}x^3 + \ldots $$
$$ = \sum_{i=0}^\infty c_i(\alpha,\beta,\gamma)x^i, $$
o\`u 
$$ c_i(\alpha,\beta,\gamma) = \frac{\alpha(\alpha + 1)\ldots (\alpha + i - 1)
\cdot \beta(\beta + 1)\ldots (\beta + i - 1)}{i!\cdot\gamma(\gamma + 1)\ldots 
(\gamma + i - 1)}, $$  
cf. [Gauss]. Il s'ensuit~: 
$$ c_i(-n/2,-n/2 + 1/2,1/2)$$
$$ = \frac{(-n/2)(-n/2 + 1)\ldots (-n/2 + i - 1)\cdot 
(-n/2 + 1/2)(-n/2 + 3/2)\ldots (-n/2 + i - 1/2)}
{i!\cdot (1/2)(1/2 + 1)\ldots (1/2 + i - 1)} $$ 
$$ = \frac{(-1)^i 2^{-i} n(n-2)\ldots (n - 2i + 2)\cdot (-1)^i
2^{-i}(n-1)(n-3)\ldots (n - 2i + 1)}{i!\cdot 2^{-i}\cdot 
1\cdot 3 \cdot 5\ldots (2i - 1)} $$ 
$$ = \frac{2^{-i}\cdot n(n-1)(n-2)\ldots (n - 2i + 1)}
{2^{-i}\cdot 2 \cdot 4\ldots 2i\cdot 1\cdot 3 \cdot 5\ldots (2i - 1)} = \binom{n}{2i}\ .$$
Donc 
$$ F(-n/2,-n/2 + 1/2,1/2,x^2) = \sum_{i=0}^{[n/2]} \binom{n}{2i}x^{2i} = 
\frac{1}{2}\{(1 + x)^n + (1 - x)^n\}\ .\eqno{(2.2.1)}$$ 
Il en d\'ecoule~:
$$ t^nF(-n/2,-n/2 + 1/2,1/2,u^2/t^2) = \frac{1}{2}\{(t + u)^n +(t - u)^n\},\eqno{(2.2.2)} $$
cf. [Gauss], no. 5, formula II.

\medskip
{\bf 2.3.} La formule (2.2.1) implique~: 
$$ E_n(x) = F(-n,-n + 1/2,1/2,- x^2/4n^2)\ .\eqno{(2.3.1)}$$

\medskip
{\bf 2.4.} Si l'on \'ecrit 
$$ E_n(x) = \sum_{k=0}^n\ e_{nk}t^{2k},\ e_{nk} := (-1)^k\binom{2n}{2k}\frac{1}{(2n)^{2k}}$$
alors 
$$e_{nk} = (-1)^k\frac{2n(2n - 1)\ldots (2n - 2k + 1)} 
{(2k)!(2n)^{2k}} = \frac{(-1)^k}{(2k)!}\cdot 
1\cdot \biggl(1 - \frac{1}{2n}\biggr) \biggl(1 - \frac{2}{2n}\biggr)\ldots 
\biggl(1 - \frac{2k-1}{2n}\biggr), $$
d'o\`u
$$ \lim_{n\ra\infty}e_{nk} = \frac{(-1)^k}{(2k)!}, $$
i.e. 
$$ \lim_{n\ra\infty} E_n(x) = \sum_{k=0}^\infty\ \frac{(-1)^k}{(2k)!}x^{2k} = \cos x, $$ 
comme il faut. En d'autres termes,
$$ \lim_{n\ra\infty}F(-n,-n + 1/2,1/2,- x^2/4n^2) = \cos x, $$
ou, comme aurait pu \'ecrire Gauss, 
$$ F(-k,k+1/2,1/2,-x^2/4k^2) = \cos x, $$
$k$ \'etant "un nombre infiniment grand" ({\it denotante $k$ numerum infinite
magnum}). En fait, Gauss \'ecrivit
$$ F(k,k',1/2,- x^2/4kk') = \cos x, $$
{\it denotante $k, k'$ numeros infinite magnos}, cf. [Gauss], 
no. 5, formula XII.    

\bigskip\bigskip

\centerline{\bf \S\  3. Asymptotiques}\label{Asymptotiques}

\bigskip 

{\bf 3.1.} On pose~: 
$$ f_{n}(x) = \sum_{k=0}^{n}\ (-1)^k\binom{2n}{2k}\frac{x^{k}}{(2n)^{2k}} = 
\sum_{k=0}^{n}\ a_k^{(n)} x^k\ .\eqno{(3.1.1)}$$
Donc 
$$ E_n(x) = f_n(x^2)\ .$$
On d\'esigne les quantit\'es $b(j)_i, r_i$, etc. 
qui correspondent au polyn\^ome $f_n$ en ajoutant l'indice $(n)$ en haut~: 
$b(j)^{(n)}_i, r^{(n)}_i$, etc. 

Donc on aura~: 
$$ c(m+1)^{(n)} = \biggl(\prod_{i=1}^{m}\ a_{n-i}^{(n)}\biggr)^2 \times  $$
$$ \times \det\left(\begin{array}{cccc} 
(n-1)\beta(1)^{(n)}_2 & (n-2)\beta(1)^{(n)}_3  & \ldots 
& (n-m)\beta(1)^{(n)}_{m+1} \\ 
(n-2)\beta(1)^{(n)}_3 & (n-2)\beta(2)^{(n)}_4  &  \ldots 
& (n-m)\beta(2)^{(n)}_{m+2}\\
 . & . & \ldots & . \\ 
(n-m)\beta(1)^{(n)}_{m+1} & (n-m)\beta(2)^{(n)}_{m+2} & 
\ldots & (n-m)\beta(m)^{(n)}_{2m} \\  
\end{array}\right)\ . $$

\medskip
{\bf 3.2.} On a~: 
$$ a_i^{(n)} = (-1)^i\binom{2n}{2i}, $$
d'o\`u
$$\begin{array}{l} r_i^{(n)} = \frac{a_i^{(n)}a_{i-2}^{(n)}}{a_{i-1}^{(n)2}} = 
\frac{[(2i-2)!]^2 [(2n-2i+2)!]^2}
{(2i)!(2n-2i)!(2i-4)!(2n-2i+4)!} \\[1em]
\ \ \ \ \ \  = \frac{(2i-2)(2i-3)}{2i (2i-1)}\cdot 
\frac{(2n-2i+1)(2n-2i+2)}{(2n-2i+3)(2n-2i+4)}\ .\end{array}$$
En rempla\c{c}ant $i$ par $n-i$, 
$$ r_{n-i}^{(n)} = 
\frac{(2i+1)(2i+2)}{(2i+3)(2i+4)}\cdot 
\frac{(2n-2i-2)(2n-2i-3)}{(2n-2i)(2n-2i-1)}\ .$$
On s'interesse aux valeurs limites~: 
$$ r_{\infty-i}^{(\infty)} := \lim_{n\rightarrow\infty} r_{n-i}^{(n)} = 
\frac{(2i+1)(2i+2)}{(2i+3)(2i+4)}\ .$$
Il s'ensuit~: 
$$\psi(\infty - i+2,\infty) := \lim_{n\rightarrow\infty} 
\psi(n - i+2,n) = \frac{1\cdot 2}{(2i-1)2i}\ ,$$
$$ \psi(\infty - i+3,\infty-1)  = \frac{3\cdot 4}{(2i-2)(2i-3)}\ ,$$
$$\psi(\infty - i+4,\infty-2)  = \frac{5\cdot 6}{(2i-4)(2i-5)}\ ,$$
etc. 

\medskip
{\bf 3.3.} Maintenant on veut calculer 
$$ \beta(j)_i^{(\infty)} := \lim_{n\rightarrow\infty} 
\beta(j)_i^{(n)}\ .$$
Il est commode de poser~: 
$$ B(j)^{\infty}_i := \beta(j)_i^{(\infty)} + j\ .$$
On a~: 
$$ B(1)_i^{(\infty)} = i \cdot \psi(\infty - i + 2,\infty) =   
\frac{1}{2i-1} $$
d'o\`u
$$ \beta(1)_i^{(\infty)} = -\frac{2(i-1)}{2i-1}\ .$$
Ensuite, 
$$ B(2)_i^{(\infty)} = i \cdot \psi(\infty - i + 2,\infty) 
\psi(\infty - i + 3,\infty - 1) + (i-2)\cdot 
\psi(\infty - i + 3,\infty - 1) $$
$$ = \psi(\infty - i + 3,\infty - 1)\cdot \biggl\{ 
B(1)_i^{(\infty)} + i - 2 \biggr\} = \frac{3\cdot 4}{(2i-2)(2i-3)}\cdot \biggl\{ 
\frac{1}{2i-1} + i - 2\biggr\} = \frac{3\cdot 2}{2i-1}, $$ 
d'o\`u 
$$ \beta(2)_i^{(\infty)} = -\frac{4(i-2)}{2i-1}\ .$$
De m\^eme, 
$$ B(3)_i^{(\infty)} = \psi(\infty - i + 4,\infty - 2)\cdot \biggl\{ 
B(2)_i^{(\infty)} + i - 4 \biggr\} 
 = \frac{5\cdot 6}{(2i-4)(2i-5)}\cdot \biggl\{ 
\frac{3\cdot 2}{2i-1} + i - 4\biggr\} = \frac{5\cdot 3}{2i-1}, $$
d'o\`u 
$$ \beta(3)_i^{(\infty)} = -\frac{6(i-3)}{2i-1}\ .$$

\medskip
{\bf 3.4.} En g\'en\'eral, la r\'ecurrence \'evidente fournit 
$$ B(j)_i^{(\infty)} = \frac{(2j-1)\cdot j}{2i-1} $$
et
$$ \beta(j)_i^{(\infty)} = - \frac{2j(i-j)}{2i-1}\ .$$

\medskip
{\bf 3.5.} On d\'efinit les nombres 
$$ \fc(m+1)^{\infty} := \det\left(\begin{array}{cccc} 
\beta(1)^{(\infty)}_2 & \beta(1)^{(\infty)}_3  & \ldots 
& \beta(1)^{(\infty)}_{m+1} \\[1em] 
\beta(1)^{(\infty)}_3 & \beta(2)^{(\infty)}_4  &  \ldots 
& \beta(2)^{(\infty)}_{m+2}\\
 . & . & \ldots & . \\ 
\beta(1)^{(\infty)}_{m+1} & \beta(2)^{(\infty)}_{m+2} & 
\ldots & \beta(m)^{(\infty)}_{2m} %\\  
\end{array}\right)\ .$$
Donc on aura~: 
$$ \biggl(\prod_{i=1}^{m}\ a_{n-i}^{(n)}\biggr)^{-2}\cdot c(m+1)^{(n)} = 
\fc(m+1)^{\infty}\cdot n^m  + O(n^{m-1})\ .$$
Les calculs pr\'ec\'edents fournissent par exemple~: 
$$ \fc(4)^{\infty} = \det\left( \begin{array}{ccc} 
- \frac{2}{3} & - \frac{4}{5} & - \frac{6}{7} \\ [0.5em] 
- \frac{4}{5} & - \frac{8}{7} & - \frac{12}{9} \\ [0.5em] 
- \frac{6}{7} & - \frac{12}{9} & - \frac{18}{11} 
\end{array}\right) = (-1)^3 \cdot 2\cdot 4\cdot 6 \cdot 
\det\left( \begin{array}{ccc} 
\frac{1}{3} & \frac{2}{5} & \frac{3}{7} \\ [0.5em] 
\frac{1}{5} & \frac{2}{7} & \frac{3}{9} \\ [0.5em] 
\frac{1}{7} & \frac{2}{9} & \frac{3}{11} 
\end{array}\right) $$
$$ = (-1)^3 \cdot 2^3 \cdot (3!)^2 \cdot 
\det\left(\begin{array}{ccc} 
\frac{1}{3} & \frac{1}{5} & \frac{1}{7} \\ [0.5em] 
\frac{1}{5} & \frac{1}{7} & \frac{1}{9} \\ [0.5em] 
\frac{1}{7} & \frac{1}{9} & \frac{1}{11} 
\end{array}\right)\ .$$
En g\'en\'eral on obtient 
$$ \fc(m+1)^{\infty} = (-1)^m \cdot 2^m \cdot (m!)^2 \cdot 
\det\left(\begin{array}{cccc} 
\frac{1}{3} & \frac{1}{5} & \ldots & \frac{1}{2m+1} \\ [0.5em] 
\frac{1}{5} & \frac{1}{7} & \ldots & \frac{1}{2m+3} \\ [0.5em] 
. & . & \ldots & . \\ [0.5em]
\frac{1}{2m+1} & \frac{1}{2m+3} & \ldots & \frac{1}{4m-1} 
\end{array}\right)\ .$$
On remarque que la derni\`ere matrice (une variante de la matrice de Hilbert) 
est du type Hankel. 

\medskip
{\bf 3.6.} Le d\'eterminant 
$$ \fC(m+1) := \det\left(\begin{array}{cccc} 
\frac{1}{3} & \frac{1}{5} & \ldots & \frac{1}{2m+1} \\ [0.5em]
\frac{1}{5} & \frac{1}{7} & \ldots & \frac{1}{2m+3} \\ [0.5em]
. & . & \ldots & . \\ [0.5em]
\frac{1}{2m+1} & \frac{1}{2m+3} & \ldots & \frac{1}{4m-1} 
\end{array}\right)$$
est un cas particulier du d\'eterminant calcul\'e par Cauchy 
(d'o\`u le caract\`ere $\fC$), cf. 
son {\it M\'emoire sur les fonctions altern\'ees et sur les sommes altern\'ees},
pp. 173 - 182 dans [Cauchy]. 

%(ou [PS], Abschitt VII, Aufgabe 3)

Rappelons que, \'etant donn\'ees deux suites $x_1, \ldots, x_m$ et 
$y_1, \ldots, y_m$, le th\'eor\`eme de Cauchy dit que 
$$\det\bigl( (x_i + y_j)^{-1}\bigr)_{i,j = 1}^m = 
\frac{\prod_{1 \leq i < j \leq m}\ (x_j - x_i)(y_j - y_i)}
{\prod_{i,j = 1}^m\ (x_i + x_j)}  $$
d'o\`u, en posant $x_i = 2i-2,\ y_i = 2i+1$, 
$$\fC(m+1) = \frac{\prod_{1 \leq i < j \leq m}\ (2j - 2i)^2}
{\prod_{i,j = 1}^m\ (2i + 2j - 1)}\ .$$

\bigskip\bigskip  

\centerline{\bf Bibliographie}\label{Bibliographie}

\bigskip

\noindent [Cauchy] Exercices d'Analyse et de Physique Math\'ematique, par le 
Baron Augustin Cauchy, Tome Deuxi\`eme, Paris, Bachelier, 1841~; {\it Oeuvres 
compl\`etes}, II-e s\'erie, tome {\bf XII}, Gauthier-Villars, MCMXVI.    

\noindent [Euler] L.Euler, De summis serierum reciprocarum ex potestatibus numerorum 
naturalium ortarum dissertatio altera in qua eaedem summationes ex fonte 
maxime diverso derivantur, {\it Miscellanea Berolinensia} {\bf 7}, 
1743, pp. 172 - 192. 

\noindent [Gauss] C.F.Gauss, Circa seriem infinitam $1 + \frac{\alpha\beta}{1 .\gamma}x + 
\frac{\alpha(\alpha + 1)\beta(\beta + 1)}{1 . 2 .\gamma(\gamma + 1)}xx + $ 
%\newline 
$\frac{\alpha(\alpha + 1)(\alpha + 2)\beta(\beta + 1)(\beta + 2)}
{1 . 2 .3 . \gamma(\gamma + 1)(\gamma + 2)}x^3 + $ etc. Pars Prior, 
{\it Commentationes societatis regiae scientarum Gottingensis recentiores}, 
Vol. {\bf II}, Gottingae MDCCCXIII. 

\noindent [Jacobi] C.G.J. Jacobi, De eliminatione variabilis e duabus
aequationibus algebraitis, {\it Crelle J. f\"ur reine und angewandte Mathematik},
{\bf 15}, 1836, ss. 101 - 124.    

\noindent [Sturm] C.-F. Sturm, M\'emoire sur la r\'esolution des \'equations num\'eriques, 
{\it M\'emoires pr\'esent\'es par divers savants \`a l'Acad\'emie Royale des 
Sciences}, Sciences math\'ematiques et physiques,  tome {\bf VI}, 1835, pp. 271 - 318.    

\end{document}